\documentclass[reqno]{amsart}

\usepackage{amssymb,latexsym, amsmath, amsxtra,amsthm}

\usepackage{graphicx}
\usepackage{array}

\numberwithin{equation}{section}

\newcommand{\cj}[1]{\overline{#1}}
\newcommand{\abs}[1]{\left| #1 \right|}
\newcommand{\cbr}[1]{\left\{ #1 \right\}}

\newcommand{\precis}{\underline}

\newcommand{\Q}{{\mathbb Q}}
\newcommand{\Z}{{\mathbb Z}}
\newcommand{\R}{{\mathbb R}}
\newcommand{\C}{{\mathbb C}}

\newcommand{\GL}{\mathrm{GL}}
\newcommand{\PSL}{\mathrm{PSL}}

\newcommand{\Sp}{{\rm Sp}}

\newcommand{\Lfunction}{\hbox{L-function}}
\newcommand{\Lfunctions}{\hbox{L-functions}}

\newtheorem{lemma}{Lemma}[section]
\newtheorem{theorem}[lemma]{Theorem}
\newtheorem{conclusion}[lemma]{Conclusion}
\newtheorem{corollary}[lemma]{Corollary}

\newtheorem*{terminology*}{Terminology}

\author{David W. Farmer}
\address{American Institute of Mathematics\\360 Portage Avenue \\Palo Alto, CA 94306-2244}
\email{farmer@aimath.org}

\author{Sally Koutsoliotas}
\address{Bucknell University\\Department of Physics and Astronomy\\Lewisburg, PA  17837} 
\email{s.koutsoliotas@bucknell.edu}

\author{Stefan Lemurell}
\address{Chalmers University of Technology and University of Gothenburg\\Department of Mathematics\\Gothenburg, Sweden}
\email{sj@chalmers.se}


\begin{document}
\title{Varieties via their L-functions}

\begin{abstract}
We describe a procedure for determining the existence, or non-existence,
of an algebraic variety of a given conductor via an analytic calculation
involving \Lfunctions.  The procedure assumes that the Hasse-Weil
\Lfunction\ of the variety satisfies its conjectured functional equation,
with no assumption of an associated automorphic object or
Galois representation.  We demonstrate the method by finding the
Hasse-Weil \Lfunctions\ of all hyperelliptic curves of conductor
less than 500.
\end{abstract}

\maketitle

\section{Introduction}
Many objects in number theory are associated to \Lfunctions.
Familiar examples include
modular forms, Galois representations, and algebraic varieties.
One might think of those as the fundamental objects, and the
associated \Lfunction\ as a shadow that serves as a placeholder
for some of the invariants of the object.  For example, the \Lfunction\ of
an elliptic curve contains information about the conductor of the curve
and the number of points on the curve $\text{mod } p$, for each prime~$p$.
But the \Lfunction\ does not know 
everything about the curve. 
For example, it does not know how many torsion points are on the
curve:  it only has information about a particular combination of the torsion order,
regulator, and a few other parameters as described in the formula of
Birch and Swinnerton-Dyer.  Of course, this apparent shortcoming is due to the
fact that the translation from an elliptic curve to its \Lfunction\ goes
via the Jacobian, and the Jacobian only sees the isogeny class of the curve.
As far as we know, no information is lost when translating from the
Jacobian/isogeny class to the \Lfunction.

While it has been natural to view \Lfunctions\ as arising
from other objects, in this paper we make the point that it
is possible to detect the \Lfunction\ 
\emph{without having access to the associated object}.
By ``detect'' we mean that one can find the Dirichlet coefficients
and functional equation parameters in a practical manner involving a computer calculation.
We have demonstrated this principle for \Lfunctions\ of
Maass forms on $\GL(3)$,
$\GL(4)$, and $\Sp(4)$, see~\cite{FKL}.  
There are many internal consistency checks
for those calculations, but in most cases there does not currently
exist an independent method to
generate the data for the underlying Maass form.
In the present paper, the \Lfunctions\ we consider are associated
to several types of arithmetic objects, as described
in Section~\ref{sec:related}.  In most cases this allows for 
independent verification of our calculations.


In the next section, we describe the \Lfunctions\  we will consider, we
list (some of) the associated objects, and we summarize the results of our calculations.
Details of our method are given in the following sections.

\section{Axioms for \Lfunctions}
This paper takes an axiomatic approach to \Lfunctions.
The most famous example of this approach is the Selberg class~\cite{Sel}.
Selberg's axioms are as general as possible: all known (or conjectured)
\Lfunctions\ satisfy much stricter conditions.
Since our methods involve computer calculations,  
we seek axioms that are
as restrictive as possible. This reduces the size of the space we must search.

A proposal for such a set of axioms is given in~\cite{FPRS}.  
We specialize to the case of
``tempered arithmetic entire \Lfunctions\ of degree~4, weight~1, and trivial central
character, which have rational integer coefficients in the arithmetic normalization.''
While that set of \Lfunctions\ may appear highly specialized, there are many objects
which give rise to such \Lfunctions, for example a hyperelliptic curve $C/\Q$ of genus~2
or an abelian surface~$A/\Q$.
See Section~\ref{sec:related} for more examples.

\subsection{The Axioms}\label{sec:axioms}
Below is a complete list of the properties we assume for the \Lfunctions\ in
our search. 

\vskip 0.1in
\noindent\textbf{\textsl{Dirichlet series:}} 
\begin{equation}
\label{dirichletseries}
L(s) = \sum_{n=1}^\infty \frac{A_n / \sqrt{n}}{n^s},
\ \ \ \ \ \ \ \ \ \ \ \ \ \ \
A_n \in \Z,
\end{equation}
with $A_n\ll n^{\frac12+\varepsilon}$.

\vskip 0.1in
\noindent\textbf{\textsl{Functional equation:}} 
\begin{equation}
\label{FE} 
\Lambda(s) := N^{s/2} \ \Gamma_\C(s+\tfrac12)^2 \, L(s),
\end{equation}
originally defined by \eqref{dirichletseries} for $s=\sigma+it$ with $\sigma>1$,
continues to an entire function which satisfies the functional equation
\begin{equation}
\Lambda(s)
 \, = \, \varepsilon \, \Lambda(1-s),
\end{equation}
where $\varepsilon = \pm 1$.
The positive integer $N$ is called the \emph{conductor} of the \Lfunction, and
$\varepsilon$ is called the \emph{sign}.

\vskip 0.1in
\noindent\textbf{\textsl{Euler Product:}} 
\begin{equation}
\label{eulerproduct} 
L(s) = \prod_p F_p(p^{-s})^{-1},
\end{equation}
with
$F_p(z)=G_p(z/\sqrt{p})$ where $G_p(z) \in \Z[z]$.
Furthermore, if $p\nmid N$ then
\begin{equation}
G_p(z) = 1 - A_p z + (A_p^2 - A_{p^2})z^2 - A_p z^3 + z^4
\end{equation} 
and the roots of $F_p(z)$ lie on $|z|=1$,
and if $p|n$ then $F_p$ has degree~$< 4$ and each  root of $F_p(z)$ lies on $|z|=p^{m/2}$
for some $m\in \{0,1,2,3\}$.

Note that the restrictions on $F_p$ imply that the Dirichlet coefficients satisfy
the bound $|A_n|/\sqrt{n} \le d_4(n)$, where $d_4(n) = \sum_{abcd=n}1$ is the
$4$-fold divisor function; in particular, $d_4(p) = 4$ if $p$ is prime.

In the functional equation we used the normalized $\Gamma$-function
\begin{equation}
\Gamma_\C(s) = 2 (2\pi)^{-s} \Gamma(s).
\end{equation}
For most of the paper we write our \Lfunctions\ as
\begin{equation}
L(s) = \sum_{n=1}^\infty \frac{b_n}{n^s},
\end{equation}
and only occasionally will we make use of the fact that
$b_n = A_n/\sqrt{n}$ for some $A_n\in \Z$.

\subsection{Arithmetic \emph{vs} analytic normalization}
In this paper we consider \Lfunctions\  in the analytic normalization,
meaning that the functional equation relates $s$ to $1-s$.
The benefit of providing a uniform approach to all \Lfunctions\  is partially offset
by the fact that, for some \Lfunctions, the arithmetic nature of the coefficients has been obscured.
Thus, it is natural to also consider the arithmetically normalized
\Lfunction, which in the case at hand is given by
\begin{equation}
L_{ar}(s) = L(s-\tfrac12) =  \sum_{n=1}^\infty \frac{A_n}{n^s},
\end{equation}
which satisfies the functional equation
\begin{equation}
\Lambda_{ar}(s) := N^{s/2} \ \Gamma_\C(s+1)^2 \, L_{ar}(s) \, = \, \varepsilon \, \Lambda_{ar}(2-s).
\end{equation}
This perspective is often preferred by people who consider
\Lfunctions\ of algebraic varieties, in particular elliptic curves and
hyperelliptic curves.  However, not all \Lfunctions\ have an arithmetic
normalization, for it is expected that all primitive  \Lfunctions\ 
$L(s) = \sum b_n n^{-s}$ fall into one of these three categories:
\begin{enumerate}

\item[(M)]{} Motivic \Lfunctions, for which there exists a number field $F/\Q$ and
an integer $w$ such that $b_n n^{w/2} \in {\mathcal O_F}$ for all $n$.

\item[(N)]{} Non-arithmetic \Lfunctions, for which a positive proportion of the Dirichlet coefficients
  are transcendental.

\item[(O)]{} Other \Lfunctions, for which there exists 
an integer $w$ such that $b_n n^{w/2}$ is an algebraic integer for all $n$,
but those integers do not all lie in the same number field.

\end{enumerate}

An example of type (O) is $L(s,\chi)L(s,E)$, where $\chi$ is a  Dirichlet character
and $E$ is an elliptic curve.
We are not aware of any examples of primitive \Lfunctions\  of type (O).
See \cite{FPRS} for more discussion.

\subsection{Related objects}\label{sec:related}
While we wish to make the point that \Lfunctions\ with the above axioms
can be thought of as existing on their own, of course
there are various objects which are associated to such \Lfunctions.
We list example objects below; we would welcome learning of additional
items to include on this list.

The following objects
are associated to
\Lfunctions\ that 
 (in some cases conjecturally)
satisfy a functional equation
of the form~\eqref{FE}:

\begin{itemize}
\item[1) ]  Siegel cusp forms of weight 2 on the paramodular group $K(N)$;

\item[2) ]  Hyperelliptic curves $C/\Q$ of genus 2 with conductor~$N$;

\item[3) ]  Abelian surfaces $A/\Q$ of conductor~$N$.
\end{itemize}

In the next three examples, $k/\Q$ is a quadratic field with ring
of integers~$\mathcal O_k$ and discriminant $D$, and
$\mathfrak{N}$ is an ideal of $\mathcal O_k$.  In each case,
the conductor of the associated \Lfunction\ is given by
$N=\operatorname{Norm}(\mathfrak{N}) D^2$.

\begin{itemize}
\item[4) ] Elliptic curves $E/k$ of conductor~$\mathfrak{N}$;

\item[5) ] Hilbert cusp forms of parallel weight $2$ 
for $\Gamma_0(\mathfrak{N}) \subset  \PSL(2,\mathcal O_k)$, where $k$ is real quadratic;

\item[6) ] Bianchi cusp forms of weight $2$ 
for $\Gamma_0(\mathfrak{N}) \subset \PSL(2,\mathcal O_k)$, where $k$ is imaginary  quadratic.
\end{itemize}

The final two cases only give rise to non-primitive \Lfunctions.
\begin{itemize}
\item[7) ] The \Lfunction\ $L(s,E_1)L(s,E_2)$, where $E_j/\Q$ is an elliptic curve
of conductor $N_j$;

\item[8) ] The \Lfunction\ $L(s,f_1)L(s,f_2)$, where $f_j \in S_2(\Gamma_0(N_j),\chi_j)$,
where either both $\chi_1$ and $\chi_2$ are the trivial character,
or $N_1=N_2$ and $\chi_2 = \overline{\chi_1}$.
\end{itemize}
In examples 7) and 8), the conductor of the \Lfunction\ is $N=N_1 N_2$.

The above examples satisfy the functional equation \eqref{FE},
but we will also require
the (arithmetically normalized) Dirichlet coefficients to be
rational integers.  This will automatically happen in cases
2), 3), 4), and 7).  In cases 1), 5), and 6), it may occur for some
cusp forms in the space but not others.  In case 8), there are
two ways to have 
rational integer coefficients in the product.
One possibility is that both $f_1$ and $f_2$ 
have rational integer coefficients, which by modularity is
just a translation of case 7) if $\chi_1$ and $\chi_2$ are trivial.
The other is that the coefficients of $f_1$ 
lie in a quadratic field and $f_2 = {f_1^\sigma}$, where $\sigma$
is the Galois conjugate.

\subsection{Modularity}
The same \Lfunction\  can arise from more than one object on the
list in Section~\ref{sec:related}.  For example, the \Lfunction\  of
a hyperelliptic curve is also the \Lfunction\ of an abelian surface,
namely, the Jacobian of the curve.  Other examples lie deeper.
For instance, by work of
Freitas, Le~Hung, and Siksek~\cite{SS}, if $k/\Q$ is real quadratic,
then every elliptic curve $E/k$ from item 7) has the same \Lfunction\ as 
that of a Hilbert modular form, item 5).  There are several other conjectural
relations, such as the paramodular conjecture~\cite{BK} which
associates certain abelian surfaces to Siegel paramodular cusp forms.

In this paper we use only the properties of the \Lfunction\ as given in Section~\ref{sec:axioms}.
In order to draw conclusions about objects related to \Lfunctions, we propose the following:
\begin{terminology*} The \emph{L-modularity conjecture} for the object $X$ is the
assertion that the \Lfunction\ of $X$ has an analytic continuation
and satisfies its conjectured functional equation.
\end{terminology*}

We note that L-modularity is a weaker condition than modularity.  
In Section~\ref{corollaries} we describe results on hyperelliptic
curves and abelian surfaces which are conditional on the L-modularity conjecture
for those objects. This is a weaker assumption than the paramodular conjecture~\cite{BK}.

\subsection{Main results}
We now state our main results, which we refer to as ``Computational Conclusion,''
or ``Conclusion'' for short.
The terminology reflects the fact that the results are based on numerical
calculations, and in principle they can be made completely rigorous.
To be made rigorous one would require a method of certifying that {a)} the computer
running the computation performed as expected, {b)} the computation was done with explicit
truncation bounds and 
error bounds (via interval arithmetic, for example),
and {c)} that the computer code was a correct
implementation of the formulas that were the theoretical basis for the computation.
Having not performed those tasks,
we do not refer to our results as a ``Theorem.''
But we emphasize that \emph{in principle} the results
can be made rigorous, which is qualitatively different than only providing 
numerical evidence.

\begin{conclusion}\label{conc:1}
The only conductors $N \le 500 $ for which there can exist an \Lfunction\ satisfying
the axioms in Section~\ref{sec:axioms} are the following:
\begin{enumerate}
\item[1)] The 70 values $N=N_1 N_2$, where $N_j$ is the conductor of
an elliptic curve~$E/\Q$.
\item[2)] The 6 values $N=M^2$, where there is an $f \in S_2(\Gamma_0(M),\chi)$ where
$\chi$ is a real nontrivial Dirichlet character and the Fourier coefficients of $f$
lie in a quadratic field.
\item[3)] The 12 values
$N=249$,    
$277$,    
$295$,    
$349$,    
$353$,    
$388$,    
$389$,    
$394$,    
$427$,    
$461$,    
$464$, and
$472$, each of which has $\varepsilon=1$.
\end{enumerate}
\end{conclusion}
Note that the values of $N$ in case~1) begin 121, 154, 165, 169, ...,
and the six values of $N$ in case~2) are $13^2$, $16^2$, $18^2$, $20^2$,
and $21^2$~(twice).

There are examples of case~2) where the character is trivial. 
According to the LMFDB~\cite{LMFDB}, the first
occurs with $N=23^2 = 529$, which is beyond the range covered by this theorem.

\begin{conclusion}\label{conc:2}
For each $N$ in Conclusion~\ref{conc:1}, there is only one choice of
Dirichlet coefficients
$b_2$, $b_3$, $b_4$, $b_5$, \ldots, $b_{100}$ for which there can exist an
\Lfunction\ satisfying
the axioms in Section~\ref{sec:axioms}, with the exception
of $N = 286$, $364$, $390$, $400$, $418$, $442$, and $480$, for which there are two choices,
and $N=441$, for which there are three choices.
\end{conclusion}

Note: the conductors in which there are multiple \Lfunctions\ in Conclusion~\ref{conc:2}
all arise from cases~1) and~2) in Conclusion~\ref{conc:1}.


\subsection{Interpretation of the results}\label{corollaries}

The results in Conclusion~\ref{conc:1} have implications for the arithmetic
objects associated to such \Lfunctions.  For example:

\begin{corollary}
Assuming L-Modularity for hyperelliptic curves, a hyperelliptic
curve $C/\Q$ must have conductor $N \ge 154$.
\end{corollary}

This improves the previous lower bound $N\ge 105$, due to Mestre~\cite{Mes},
under the same assumptions.

\begin{proof} The only value $N< 154$ permitted by Conclusion~\ref{conc:1}
is $N=121$.  But Schoof proved~\cite{Scho} that the only curve$/\Q$ which has good reduction outside
of $p=11$ is an elliptic curve of conductor~11.
\end{proof}

There are hyperelliptic curves with conductor 169,
such as $y^2 = x^6 + 4x^5+6x^4+2x^3+x^2+2x+1$ on Stoll's list~\cite{Stol}.
In order to prove that the smallest conductor is $N=169$, one must
show that 154 and 165, which do arise as the product of two elliptic
curves, do not also arise from hyperelliptic curves.
Conclusion~\ref{conc:2}~states that for each of those conductors there is a unique choice
of the first 100 Dirichlet coefficients, and those agree with the
coefficients of $L(s,E_{11})L(s,E_{14})$ and $L(s,E_{11})L(s,E_{15})$, respectively.
It is possible to determine more coefficients, which could allow one to apply
the Faltings-Serre method to prove that those are the only \Lfunctions\  with those
conductors arising from an algebraic variety.
It would still remain to prove that there is no hyperelliptic curve with
such an \Lfunction.

Other results along these lines can be obtained.  For example, Conclusion~\ref{conc:1}
can be used to give a new proof of the result~\cite{cremona} that $25$ is 
the smallest norm conductor of an elliptic curve$/\Q(i)$.

We can also use Conclusion~\ref{conc:1} to resolve some of the cases in Table~1
in the work of Brumer and Kramer~\cite{BK}.  Brumer and Kramer
seek to produce a provably complete list of all conductors of abelian surfaces.
For every
odd integer $N<1000$ 
they either have an example of an abelian surface with such a conductor,
or 
they prove there is
no such surface.  However, there are 50 cases of odd $N<1000$
which they were unable to resolve,
the first few of which are $N=415$, 417, 531, and 535.

An immediate consequence of Conclusion~\ref{conc:1} is:

\begin{corollary} Assuming L-Modularity for abelian surfaces, there
is no abelian surface with conductor $N=415$ or $417$.
\end{corollary}

\subsection{Comparison with other tables of arithmetic objects}\label{sec:efforts}

We briefly discuss the relationship between our calculations and other
work concerning the objects listed in Section~\ref{sec:related}.

Brumer and Kramer's~\cite{BK}  formulation of the
Paramodular Conjecture was greatly aided by the calculations of
Poor and Yuen~\cite{PY} of paramodular cusp forms of weight~2 and prime level.
Poor and Yuen (personal communication) have extended their calculations
to all levels $N<1000$.  Mainly due to a lack of dimension formulas, their 
results are rigorous only in a few cases.  However, their heuristic results are consistent with
our calculations, and also with the tabulation of hyperelliptic curves
that we describe next.

Booker, Sijsling, Sutherland, Voight, and Yasaki~\cite{BSSVY} tabulated 
hyperelliptic curves$/\Q$ of genus~2.  They found hyperelliptic curves
with each of the 12 conductors listed in item~3) of Conclusion~\ref{conc:1},
and they have not found any conductors $N<500$ which are not listed in Conclusion~\ref{conc:1}.
(There are conductors arising from abelian surfaces that are not isogenous
to the Jacobian of a hyperelliptic curve, but the smallest known has conductor $N=561$.)

With the exception of the first elliptic curve$/\Q(i)$ 
(whose L-function has conductor $N=25\cdot 4^2 = 400$), 
the objects listed in Section~\ref{sec:related} have conductors
outside the range of our calculations.  

An interesting case is conductor $N=550$, for which Poor and Yuen~\cite{PY2}
have found a paramodular form, but for which the Abelian surface
predicted by the Paramodular Conjecture has not yet been found.
In Section~\ref{sec:550} we provide data for the \Lfunction, which
may be of some use in the search for the conjectured Abelian surface.

\section{A sketch of the method}\label{sec:sketch}
We briefly describe the methods used to prove Conclusion~\ref{conc:1}.
The details comprise the rest of the paper.

For a given conductor $N$ and sign of the functional equation~$\varepsilon$,
we wish to determine whether there is an \Lfunction\ which satisfies 
a functional equation with those parameters as well as the other conditions
described
in Section~\ref{sec:axioms}.  At this point, the only information we don't know about
the \Lfunctions\  are the Dirichlet coefficients.

We treat the Dirichlet coefficients as unknowns, and form a system of
equations involving those coefficients.  The equations arise by evaluating
the \Lfunction\ in several ways using the approximate functional equation,
described in Section~\ref{sec:appfe}.
The method used to form equations is given in Section~\ref{sec:makeeqns}.

The idea of using the approximate functional equation to generate
equations with the Dirichlet coefficients as unknowns has appeared
previously.  Booker~\cite{Book} used this to find missing coefficients
when several other coefficients were already known, and Bian~\cite{Bian}
used the functional equation for twists to generate a large system
of linear equations, from which all the coefficients could be calculated.

Initially the equations are linear in the Dirichlet coefficients.
We construct a nonlinear system
by using the fact that the coefficients are multiplicative and 
satisfy a recursion for the prime powers arising from the shape of the local
factors in the Euler product.
It is possible to directly solve the nonlinear system, as we did in our previous
work on Maass form \Lfunctions~\cite{FKL}.  However, in the case at hand, we
can exploit the fact that the arithmetically normalized coefficients are
rational integers, which
implies that there are only finitely many
choices for the local factor at each~$p$.  

We compiled a complete list of possible local factors for each small prime.
This involved using the bound on the roots to obtain bounds on the coefficients,
looping over all polynomials with integer coefficients in the given range,
and then checking which had roots with the appropriate absolute value.
In Table~\ref{table:numlocalfactors}
we show the number of possibilities for small~$p$.

\begin{table}[htb]
\centering
\begin{tabular}{c c c}
\hline \hline
prime, $p$ & \ \ \ good \ \ \ & \ \ \ bad \ \ \ \\
\hline
2            & 35 & 26     \\
3            & 63 & 32  \\
5            & 129 & 38  \\
7            & 207 & 44  \\
\hline
\\
\end{tabular}
\caption{Number of possible local factors $F_p$}
\label{table:numlocalfactors}
\end{table}

It is noteworthy that in our approach, the ``bad'' primes pose no extra
difficulty.  In fact, it is easier for our method to handle the bad primes because,
as shown in Table~\ref{table:numlocalfactors}, there are fewer possibilities to
consider.

Since the system has infinitely many unknowns, we truncate it and track the
error due to the omitted terms.  
We solve the truncated system of equations by testing all possible values of
the remaining coefficients (see Section~\ref{sec:solving}). This is organized as a breadth-first search of a
tree, as described in Section~\ref{sec:tree}. 

The final ingredient is optimizing the system of equations as the search
progresses to successive depths of the tree: we use only one equation
at each conductor; see Section~\ref{sec:tree}.

\section{The approximate functional equation}\label{sec:appfe}

The  approximate functional equation
is the primary tool used to evaluate $L$-functions.
It involves a test function which
can be chosen with some freedom.  This will play a key role
in our calculations.

The material in this section is taken from Section~3.2 of~Rubinstein~\cite{Rub}.

Let
\begin{equation}
   L(s) = \sum_{n=1}^{\infty} \frac{b_n}{n^s}
\end{equation}
be a Dirichlet series that converges absolutely in a half plane, $\Re(s) > \sigma_1$.

Let
\begin{equation}
    \label{eqn:lambda}
    \Lambda(s) = Q^s
                 \left( \prod_{j=1}^a \Gamma(\kappa_j s + \lambda_j) \right)
                 L(s),
\end{equation}
with $Q,\kappa_j \in {\mathbb{R}}^+$, $\Re(\lambda_j) \geq 0$,
and assume that:
\begin{enumerate}
    \item  $\Lambda(s)$ has a meromorphic continuation to all of ${\mathbb{C}}$ with
           simple poles at $s_1,\ldots, s_\ell$ and corresponding
           residues $r_1,\ldots, r_\ell$.
    \item $\Lambda(s) = \varepsilon \cj{\Lambda(1-\cj{s})}$ for some
          $\varepsilon \in {\mathbb{C}}$, $|\varepsilon|=1$.
    \item For any $\sigma_2 \leq \sigma_3$, $L(\sigma +i t) = O(\exp{t^A})$ for some $A>0$,
          as $\abs{t} \to \infty$, $\sigma_2 \leq \sigma \leq \sigma_3$, with $A$ and the constant in
          the `Oh' notation depending on $\sigma_2$ and $\sigma_3$. \label{page:condition 3}
\end{enumerate}

Note that~\eqref{eqn:lambda} expresses the functional equation in
more general terms than~\eqref{FE}, but it is a simple matter
to unfold the definition of~$\Gamma_\R$ and~$\Gamma_\C$.

To obtain a smoothed approximate functional equation with desirable
properties, Rubinstein \cite{Rub} introduces an auxiliary function.
Let $g: \C \to \C$ be an entire function that, for fixed $s$, satisfies
\begin{equation}\label{eqn:gbound}
    \abs{\Lambda(z+s) g(z+s) z^{-1}} \to 0
\end{equation}
as $\abs{\Im{z}} \to \infty$, in vertical strips,
$-x_0 \leq \Re{z} \leq x_0$. The smoothed approximate functional
equation has the following form.
\begin{theorem}\label{thm:formula}
For $s \notin \cbr{s_1,\ldots, s_\ell}$, and $L(s)$, $g(s)$ as above,
\begin{equation}\label{eqn:formula}
  \Lambda(s) g(s) =
         \sum_{k=1}^{\ell} \frac{r_k g(s_k)}{s-s_k}
         + Q^s \sum_{n=1}^{\infty} \frac{b_n}{n^s} f_1(s,n) 
         + \varepsilon Q^{1-s} \sum_{n=1}^{\infty} \frac{\cj{b_n}}{n^{1-s}} f_2(1-s,n)
\end{equation}
where
\begin{align}\label{eqn:mellin}
   f_1(s,n) &:= \frac{1}{2\pi i}
                   \int_{\nu - i \infty}^{\nu + i \infty}
                    \prod_{j=1}^a \Gamma(\kappa_j (z+s) + \lambda_j)
                    z^{-1}
                    g(s+z)
                    (Q/n)^z
                    dz \notag \\
  f_2(1-s,n) &:= \frac{1}{2\pi i}
                   \int_{\nu - i \infty}^{\nu + i \infty}
                    \prod_{j=1}^a \Gamma(\kappa_j (z+1-s) + \cj{\lambda_j})
                    z^{-1}
                    g(s-z)
                    (Q/n)^z
                    dz
\end{align}
with $\nu > \max \cbr{0,-\Re(\lambda_1/\kappa_1+s),\ldots,-\Re(\lambda_a/\kappa_a+s)}$.
\end{theorem}

In our examples, $L(s)$ continues to an entire function, so the first
sum in \eqref{eqn:formula} does not appear.  For fixed
$Q,\kappa,\lambda,\varepsilon$, and sequence~$b_n$, and $g(s)$ as described
below, the right side of \eqref{eqn:formula}
can be evaluated to high precision.

A reasonable choice for the weight function is
\begin{equation}\label{eqn:test}
g(s)=e^{i b s + c s^2},
\end{equation}
which by Stirling's formula satisfies \eqref{eqn:gbound}
if  $c>0$, or if $c=0$ and $|b| < \pi d/4$, where $d$ is the degree of
the $L$-function.
Rubinstein~\cite{Rub} uses such a weight function with $b$ chosen
to balance the size of the terms in the approximate functional equation,
minimizing the loss in precision in the calculation.
In this paper, we exploit the fact that there are many choices of weight function,
and so there are many ways to evaluate the $L$-function.
We combine those calculations to extract as much information
as possible from the known Dirichlet coefficients.  This idea is
described in the next section.

In our discussions below, we find it more convenient to
use the Hardy $Z$-function instead of the
$L$-function itself.  The function $Z$ associated to an $L$-function
$L$ is defined by the properties: $Z(\frac12 + it)$ is a smooth function
which is real if $t$ is
real, and $|Z(\frac12 + it)| = |L(\frac12 + it)|$.

\section{A system of equations in the Dirichlet coefficients}\label{sec:makeeqns}
Once we choose a sign $\varepsilon$ and a conductor $N$,
we know everything about the \Lfunction\ except for the integers $A_n$
determining the Dirichlet coefficients.  Note that, by the
conditions on the local factors in the Euler product, the $A_n$
are determined by $A_p$ and $A_{p^2}$ for $p\nmid N$, and
 $A_p$, $A_{p^2}$, and  $A_{p^3}$ for $p|N$.

We will use the approximate functional equation~\eqref{eqn:formula} to
create a system of equations whose solution(s) will be the $A_n$.
We create the equations by using the flexibility of choosing the weight
function~$g(s)$.
The following example illustrates the idea.

In our example we let $N=211$ and $\varepsilon=1$.


In~\eqref{eqn:formula}, choose $s=\frac12 + 2 i$ and $g(z)=1$ to obtain
\begin{align}\label{eqn:choice1}
Z(\tfrac12 + 2i) =\mathstrut 
\input{Lbeta1.out}.
\end{align}

A comment on precision: numerical calculations were done to 50 decimal digits
of precision, and all numerical values shown are truncations of the actual
values.  The notation $3.308\precis{5}$ means $3.3085...$, where the
remaining digits have been truncated.

Note that since we know the exact form of the functional equation, and we have
chosen $s$ and $g(z)$, the only ``unknowns'' in the approximate functional
equations are the Dirichlet coefficients $b_n=A_n/\sqrt{n}$.

Now, choose $s=\frac12 + 2 i$ and $g(z)=e^z$ to obtain
\begin{align}\label{eqn:choice2}
Z(\tfrac12 + 2i) =\mathstrut
\input{Lbeta2.out}.
\end{align}

Since the left sides of \eqref{eqn:choice1} and \eqref{eqn:choice2} are
equal, subtracting gives the following equation, satisfied by the coefficients of
\emph{any} \Lfunction\ of the chosen form:
\begin{align}\label{eqn:subtracted}
0 =\mathstrut
\input{eqn1.out}.
\end{align}

We have the freedom to create more equations by choosing different pairs of test functions,
or by evaluating $Z(s)$ at a different point.  We establish Conclusion~\ref{conc:1} by
creating a system of such equations for each conductor~$N$ and sign~$\varepsilon$ in
the functional equations, and then solving the system subject to the constraints
imposed by the Euler product.  In the next section, we describe some of the important
properties of the system of equations, and then in Section~\ref{sec:solving} we
describe 
the process of solving the system in more detail.

\subsection{Decreasing contributions of the coefficients}
We must truncate the system of equations so that we deal only with
finitely many unknowns.  This requires bounding the contribution of those
entries which have been eliminated.  There are two ingredients to that bound:
an estimate for the size of the Dirichlet coefficients and an estimate for the
size of the terms appearing in the equations.

\begin{figure}[htp]
\begin{center}
\scalebox{0.8}[0.8]{\includegraphics{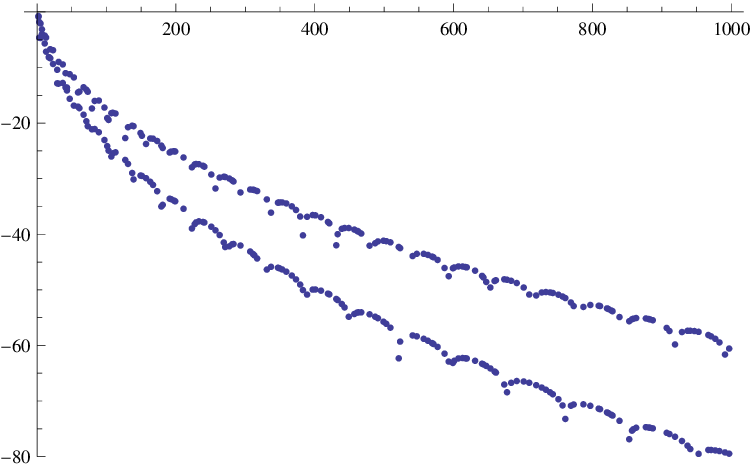}}
\caption{\sf
The ($\log$ of the) contributions of the $p^{\text{th}}$ coefficients, for $p<1000$,
 for two
example equations in the case of functional equation~\eqref{FE} with
$N=211$.
} \label{fig:plot1000}
\end{center}
\end{figure}

The Dirichlet coefficients are easily estimated using the properties we
have assumed for the Euler product.  
Our numerical examples in the previous section illustrate that
the contributions of the coefficients are decreasing fairly rapidly.  In fact, 
the size of the term multiplying $b_n$ in an equation such as
\eqref{eqn:subtracted} is approximately $\exp(C \sqrt{n/N^{1/2}})$
where $C<0$ depends on the specific $\Gamma$-factors and the test function~$g$.
This follows from the inverse Mellin transform
\begin{equation}\label{eqn:invmellin}
\frac{1}{2\pi i}
\int_{(c)} \Gamma(k z) X^{z} \, dz = e^{X^{1/k}}
\end{equation}
and the fact that $f_j(s,n)$ in \eqref{eqn:mellin}, in the case of a degree~$d$
\Lfunction,
is a perturbation of \eqref{eqn:invmellin} with $k=d/2$.  
As illustrated in Figure~\ref{fig:plot1000}, this approximation can be
fairly good even when $n$ is small, recalling that in our case $d=4$.

Proving precise estimates that are useful for small $n$ can be difficult:
see~\cite{molin} for details.  

%

\section{Solving the system}\label{sec:solving}

We continue with \eqref{eqn:subtracted}, which must be satisfied by any
\Lfunction\ with conductor $N=211$ and sign $\varepsilon=1$.

Since $2\nmid 211$, there are 35 possible
local factors at $p=2$ in the Euler product, as indicated in Table~\ref{table:numlocalfactors}.
Here are two of them:
\begin{align}
L_{2,9}(z) =\mathstrut& \left(1+\frac{1}{\sqrt{2}} \,z+ \frac{1}{\sqrt{2}}\,z^3+ z^4 \right)^{-1} \cr
     =\mathstrut& 1- \frac{1}{\sqrt{2}}\,z + \frac{1}{{2}}\,z^2  - \frac{3}{\sqrt{2}}\, z^3 
           + \frac{1}{{4}}\,z^4  + \frac{1}{4 \sqrt{2}}\, z^5 + \cdots\\
L_{2,29}(z) =\mathstrut& \left( 1-{\sqrt{2}} \,z+ z^2-{\sqrt{2}}\,z^3+ z^4\right)^{-1}\cr
     =\mathstrut&  1+ {\sqrt{2}}\,z + \,z^2  + {\sqrt{2}}\, z^3
           + {{2}}\,z^4  + {\sqrt{2}}\, z^5 + \cdots .
\end{align}
The indices 9 and 29 indicate where those polynomials fall
among the list of 35 possible local factors, when ordered lexicographically.
Note that $b_{2^m}$ is the coefficient of $z^m$.

When we substitute the two examples above into~\eqref{eqn:subtracted},
use the multiplicative relations among the coefficients, and gather
terms, we find
\begin{align}
L_{2,9}:\phantom{xxxxx}0 =\mathstrut 
\input{subs1.out}
\label{eqn:subs1}\\
L_{2,29}:\phantom{xxxxx}0 =\mathstrut
\input{subs2.out} .
\label{eqn:subs2}\nonumber
\end{align}
If we now use the bounds $|b_p|\le 4$ and $|b_{p^2}|\le 10$,
and write ``$\pm X$'' to indicate a real number in the interval $[-X,X]$,
we find
\begin{align}
L_{2,9}:\phantom{xxxxx}0 =\mathstrut& 1.106\precis{0} \pm 2.067\\
L_{2,29}:\phantom{xxxxx}0 =\mathstrut& 3.045\precis{5} \pm 2.525 \, .
\end{align}
The first equation is true, and the second equation is false.
Thus, the local factor $L_{2,29}$ cannot occur in any \Lfunction\ satisfying
our assumptions.  By creating several equations along the lines described
above, one can prove that 15 out of the 35 possible local factors
at $p=2$ \emph{cannot} occur in an \Lfunction\ satisfying our assumptions.

What can we do about the 20 local factors at $p=2$ that appear to be possible?  We move on to
the possible local factors at $p=3$.  That is, for each of the 20 remaining factors
at $p=2$, 
we consider each of the 63 possible factors at $p=3$.  To give one more numerical example,
if we use the $25^{th}$ choice for the local factor at $p=3$ in \eqref{eqn:subs1} we find:
\begin{align}
L_{2,9}\text{ and }L_{3,25}:\phantom{xxxxx}0 =\mathstrut&
0.87817\precis{6}  + 0.096\precis{4} \, b_{5} + 0.025\precis{4} \, b_{7} + 0.0018\precis{1} \, b_{11}\cr 
&+ 0.0004\precis{6} \, b_{13} 
 + 0.0000\precis{2} \, b_{17} + 3.7\precis{4}\times 10^{-6} \, b_{19} \cr  
&\cdots .
\end{align}
Since $|b_p|\le 4$, we see that the above equation cannot be satisfied by the coefficients of
any \Lfunction.  This rules out the combination $L_{2,9} L_{3,25}$ in the Euler product.

As we proceed in this manner, one of two things can happen.  We may find that \emph{every}
choice of coefficients for a given $p$ leads to an equation which cannot be satisfied.
That proves there does not exist an \Lfunction\ with the given functional equation, and therefore
there does not exist any object which would have an \Lfunction\ with such a functional equation.

Or, we could find that, after searching through the first 100 coefficients, there are
combinations which are consistent with the equations we form.  This proves that only
those combinations of coefficients are possible for any \Lfunction\ with the given functional
equation.  However, it does not prove that such an \Lfunction\ exists.

We now give more details about how we organize the search through the coefficients.

\section{Searching the tree}\label{sec:tree}
The procedure described in Section~\ref{sec:solving} 
involves a breadth-first search of a tree, where we prune
nodes which have been found to be inconsistent with a possible solution
to the system of equations.

We use the entire local factor at the primes $p=2$ and $3$. For
larger $p$, we search through the possible $b_p$ and $b_{p^2}$ separately.
This is helpful because, for example, when considering $b_5$, the value
of $b_{25}$ proportionally contributes very little to the sum: see the
plots in Figure~\ref{fig:plot1000}.  Thus, we need only consider the
17 possible values of $b_5$, instead of the 129 possible values of $F_5$.

Since the running time of the tree search is proportional to the number
of equations used when testing candidate coefficients, it is desirable
to use few equations.  That is, having a small number of equations is preferred
provided that the equations are
useful for testing the candidate coefficients.

\begin{figure}[htp]
\begin{center}
\scalebox{0.8}[0.8]{\includegraphics{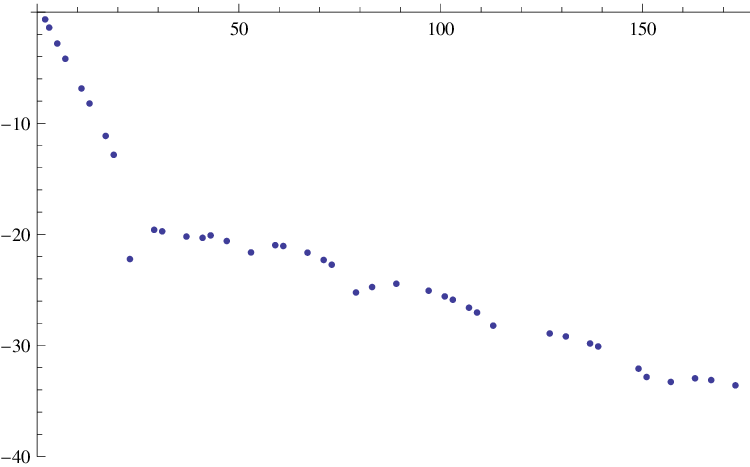}}
\caption{\sf
The ($\log$ of the) contributions of the $p^{\text{th}}$ coefficients
for one equation made from
a linear combination of 11 of our basic equations. 
} \label{fig:plotavg}
\end{center}
\end{figure}

The logical extreme, which was previously employed in~\cite{FR},
 is to construct a single equation which has been optimized
to minimize the contributions of those coefficients that contribute
to the error term.
This is easily accomplished (using least squares, for example)
by taking a linear combination
of the available equations.  Figure~\ref{fig:plotavg} illustrates an
example where we have chosen to minimize the contributions of the coefficients
$b_p$ for $p\ge23$.  As we are testing possible values of $b_{19}$ using that
equation, the error term from the larger coefficients will be small.
So one should expect that at most one
of those possible choices of $b_{19}$ will lead to a consistent equation.
See~\cite{FR} for more discussion on this approach.

To illustrate the effectiveness of the method, the following is a summary of  two sample runs.
Each list contains the number of viable possibilities at each depth of searching
the tree.  For example, in the case of $N=209$ and $\varepsilon = 1$, when considering only $p=2$,
there were 20 possible values (out of 35) that led to a consistent equation.
Then after testing every possible value with $p=3$ for each of those 20 possibilities,
a total of 87 combinations remained.  Then there were 44 possible combinations for
$p=2$, 3, and 5 that survived, and so on. 
\begin{align}
N=209,\  \varepsilon=1 :
\ \ \ \ \ \ \ \ \ \ \ \ \ \ 
&\{20, 87, 44, 4, 1, 1, 1, 1, 1, 1, 1, 1, 1, 2, 1, 1, 1, 1, 1, 1, 1, 1\}\label{conductor209}\\
N=211,\  \varepsilon=1 :
\ \ \ \ \ \ \ \ \ \ \ \ \ \ 
&\{20, 91, 39, 1, 1, 0\}\label{eliminate211}
\end{align}

From \eqref{conductor209} it appears that there is an \Lfunction\ with $N=209$ and $\varepsilon=1$.
Indeed, $L(s,E_{11})L(s,E_{19})$ is such an \Lfunction.
From \eqref{conductor209} it appears that there is exactly one \Lfunction,
but as mentioned in Section~\ref{corollaries}, the most we can say is that there is a
unique choice of the first 100 coefficients.  Current methods do not make it possible
to prove there are not several \Lfunctions\  which all begin with those same 100 coefficients.

From~\eqref{eliminate211} we see that after considering all possible values of $b_{11}$,
we have proven that there is no \Lfunction\ with $N=211$ and $\varepsilon=1$.

\section{Conductor 550}\label{sec:550}

Poor and Yuen~\cite{PY2} have identified a paramodular form of
level 550, which by the Paramodular Conjecture~\cite{BK} is associated
to an Abelian surface of conductor 550.  Such an Abelian surface
has yet to be identified, so to assist in its search, we provide
some data about the \Lfunction\  of conductor 550 which arose in our search.

Table~\ref{tab:550} provides the first few local factors
or Dirichlet coefficients, using the notation in Section~\ref{sec:axioms}.
The local factor at $p=3$ agrees with that found by Poor and Yuen.

\begin{table}[htb]
\centering
\begin{tabular}{l l || l l r }
\hline 
 $p$ & $G_p(T)$ & & $p$ & $A_p$ \\
\hline
2            & $(1 + T) (1 + 2 T^2)$   & \phantom{x} & 17 & $-3$  \\
3            & $1 - T^2 + 9 T^4$  & \phantom{x}&19  & $1$  \\
5            & $1 + 3 T + 5 T^2$  & \phantom{x}&23  & $-3$ \\
7            & $1 + 4 T^2 + 49 T^4$  & \phantom{x}&29  & $0$ \\
11            & $(1 + T) (1 - 3 T + 11 T^2)$  & \phantom{x}&31  & $-5$ \\
13            & $1 - 8 T^2 + 169 T^4$  & \phantom{x}&37  & $3$ \\
\end{tabular}
\caption{Local factors and Dirichlet coefficients of the primitive L-function of degree 4, weight 1,  and conductor $N=550$}
\label{tab:550}
\end{table}

\bibliographystyle{plain}

\begin{thebibliography}{99}

\bibitem{Book}  A. Booker, Numerical tests of modularity. J. Ramanujan Math. Soc. 20 (2005), no. 4, 283–339.

\bibitem{BSSVY} A. Booker, J. Sijsling, A. Sutherland, John Voight, and Dan Yasaki,
A database of genus 2 curves over the rational numbers, Algorithmic Number Theory 12th International Symposium (ANTS XII), LMS Journal of Computation and Mathematics 19 (2016), 235-254.

\bibitem{Bian} C. Bian, Computing GL(3) automorphic forms. 
Bull. Lond. Math. Soc. 42 (2010), no. 5, 827–842.

\bibitem{BK} A. Brumer and K. Kramer, 
Paramodular abelian varieties of odd conductor. 
Trans. Amer. Math. Soc. 366 (2014), no. 5, 2463-2516. 

\bibitem{cremona}  J. E. Cremona, Hyperbolic tessellations, modular symbols, and elliptic curves over complex quadratic fields. Compositio Math. 51 (1984), no. 3, 275-324.

\bibitem{FKL} D. Farmer, S. Koutsoliotas, and S. Lemurell, Maass forms on GL(3) and GL(4),
Int Math Res Notices 2014 (22): 6276-6301. doi: 10.1093/imrn/rnt145 
arXiv:1212.4545 

%

\bibitem{FPRS} D. Farmer, A. Pitale, N. Ryan, and R. Schmidt, Analytic L-functions: Definitions, Theorems, and Connections, preprint, arXiv:1711.10375.

\bibitem{FR} D. Farmer and N. Ryan, Evaluating L-functions with few known coefficients,
 LMS J. Comput. Math. 17 (2014), no. 1, 245-258. arXiv:1211.4181

\bibitem{LMFDB}  The LMFDB Collaboration, \emph{The L-functions and Modular Forms Database}, 
http://www.LMFDB.org [Online: accessed November 20, 2014].

\bibitem{SS} N. Freitas, B.  Le Hung, and S. Siksek 
    Elliptic Curves over Real Quadratic Fields are Modular
 arXiv:1310.7088

\bibitem{Mes} J.-F. Mestre, 
Formules explicites et minorations de conducteurs de vari\'et\'es alg\'ebriques.
Compositio Math. 58 (1986), no. 2, 209-232. 

\bibitem{molin} P. Molin, Int\'egration num\'erique et calculs de fonctions {L},
  PhD Thesis, Institut de Math\'ematiques de Bordeaux, 2010.

\bibitem{PY} C. Poor and D. S. Yuen, Paramodular Cusp Forms, arXiv:0912.0049

\bibitem{PY2} C. Poor and D. S. Yuen, personal communication.

\bibitem{Rub} M. Rubinstein, Computational methods and experiments in
analytic number theory, in, ``Recent perspectives in random matrix
theory and number theory'', F.Mezzadri and N.C.Snaith, Eds,
LMS 2005.

\bibitem{Scho} Schoof, Ren\'e Abelian varieties over Q with bad reduction in one prime only. Compos. Math. 141 (2005), no. 4, 847-868. 

\bibitem{Sel} A.  Selberg, {\it  Old and new results and conjectures about
a class of Dirichlet series}, Proceedings of the Amalfi Conference
on Analytic Number Theory 1989,
Univ. Salerno, Salerno (1992), 367--385.


\bibitem{Stol} http://www.mathe2.uni-bayreuth.de/stoll/

\end{thebibliography}

\end{document}